\documentclass[titlepage,twoside,12pt]{article}
\usepackage{amssymb}
\usepackage{amsfonts}
\newcommand{\D}{\displaystyle}

\begin{document}

{\LARGE \bf Genuine Lie semigroups} \\

{\LARGE \bf and semi-symmetries of PDEs} \\ \\

{\bf Elem\'{e}r E ~Rosinger} \\ \\
{\small \it Department of Mathematics \\ and Applied Mathematics} \\
{\small \it University of Pretoria} \\
{\small \it Pretoria} \\
{\small \it 0002 South Africa} \\
{\small \it eerosinger@hotmail.com} \\ \\

{\bf Abstract} \\

Any Lie group $G$ acting on a Euclidean nonvoid open subset $M$ can be seen as a subgroup of
the smooth diffeomorphisms ${\cal D}iff^\infty ( M, M )$ of $M$ into itself. Thus actions by
such Lie groups $G$ correspond to smooth {\it coordinate transforms} on $M$ which, in
particular, have smooth inverses. \\
In Rosinger [1, chap. 13], the study of Lie {\it semigroups} $G$ in the {\it vastly larger}
semigroup ${\cal C}^\infty ( M, M )$ of smooth maps of $M$ into itself was initiated. Such
semigroups were named {\it genuine Lie semigroups}, or in short, GLS, since they are {\it no
longer} contained in ${\cal D}iff^\infty ( M, M )$, thus they correspond to smooth coordinate
transforms which {\it need not} have smooth inverses. \\
Genuine Lie semigroups, or GLS, have a major interest since they still can transform solutions
of linear or nonlinear PDEs into other solutions of the respective equations, thus leading to
the vastly larger class of {\it semi-symmetries} of such equations. \\
Certain Lie semigroups have earlier been studied in the literature, Hilgert, et.al. However,
such semigroups have always been contained in suitable Lie groups, thus they have been
contained in ${\cal D}iff^\infty ( M, M )$ as sub-semigroups. In particular, the coordinate
transforms defined by them were always invertible, unlike in the case of genuine Lie
semigroups, or GLS, studied here. \\ \\

{\bf 1. Motivation} \\

Given an open subset $M \subseteq \mathbb{R}^l$, it is obvious that a Lie group $G$ acting on
$M$, that is \\

(1.1) $~~~ G \times M \ni ( g, x ) ~\longmapsto~ g x \in M $ \\

can be identified with a subgroup of all the smooth diffeormorphism of $M$. Namely, we have
the {\it injective group homomorphism} \\

(1.2) $~~~ G \ni g ~\longmapsto~ f_g \in {\cal D}iff^\infty ( M, M ) $ \\

where $f_g$ is defined by \\

(1.2$^*$)~~~ $ M \ni x ~\longmapsto~ f_g ( x) = g x \in M $ \\

Here the noncommutative group structure on ${\cal D}iff^\infty ( M, M )$ is defined by the
usual composition of mappings, and thus the neutral element is $e = id_M$, that is, the
identity mapping of $M$ onto itself. \\

In this way, in terms of the Euclidean domain $M$, the group homomorphism (1.2) is but a group
of smooth coordinate transforms which have smooth inverses. \\

Given a linear or nonlinear PDE \\

(1.3) $~~~ T ( x, D )~ U ( x ) ~=~ 0,~~~ x \in \Omega $ \\

where $\Omega$ is an open subset in $\mathbb{R}^n$, one of the major interests in Lie groups -
according to Lie's original aim - is in the study of the symmetries of solutions $U : \Omega
\longrightarrow \mathbb{R}$ of (1.3), which therefore, transform them into other solutions of
(1.3). This can be done as follows. One takes $M = \Omega \times \mathbb{R}$ and finds the
corresponding Lie group actions (1.1) which, when extended to the solutions $U \in {\cal
C}^\infty ( \Omega, \mathbb{R} )$ of (1.3), will transform them into solutions of the same
equation,  Appendix, or Olver, Rosinger [1]. \\

We present here a significant {\it extension} of this classical symmetry method. For that
purpose we note the following {\it four} facts :

\begin{itemize}

\item
In view of (1.2), for every Lie group $G$ acting on any nonvoid open subset $M \subseteq
\mathbb{R}^l$, we have the inclusions

\end{itemize}

(1.4) $~~~ G ~\subseteq~ {\cal D}iff^\infty ( M, M ) ~\subset~ {\cal C}^\infty ( M, M ) $

\begin{itemize}

\item Even in the simplest one-dimensional case when $M = \mathbb{R}$, the set ${\cal
C}^\infty ( M, M )$ is {\it considerably larger} than ${\cal D}iff^\infty ( M, M )$. Indeed,
only those functions $f : M \longrightarrow M$ in ${\cal C}^\infty ( M, M )$ belong to
${\cal D}iff^\infty ( M, M )$ which are bijective and have their derivative either everywhere
strictly positive, or everywhere strictly negative on $M$.

\item
With respect to the usual composition of mappings, ${\cal C}^\infty ( M, M )$ is a
noncommutative {\it semigroup} with identity, which is {\it not} a group since it contains a
vast amount of non-invertible mappings, while on the other hand, ${\cal D}iff^\infty ( M, M )$
is a noncommutative subgroup of it, and as such, the largest one.

\item
As seen in the sequel, the property that solutions $U \in {\cal C}^\infty ( \Omega,
\mathbb{R} )$ of PDEs in (1.3) are transformed by suitable Lie group actions $g \in G$ - thus
in view of (1.4), by diffeomorphism in ${\cal D}iff^\infty ( M, M )$ - into other solutions of
(1.3), is by {\it no means} restricted to diffeomorphisms alone, but it can also be valid for
certain smooth and {\it not} necessarily invertible transformations which are elements in the
considerably larger set ${\cal C}^\infty ( M, M )$.

\end{itemize}

{~}

{\bf 2. Genuine Lie Semigroups} \\

Our main aim is to {\it extend} the usual Lie group actions (1.1), (1.2), that is

$$ G \times M \ni ( g, x ) ~\longmapsto~ g x \in M,~~~~~
                G \ni g ~\longmapsto~ f_g \in {\cal D}iff^\infty ( M, M ) $$

\medskip
in such a way that we are {\it no longer} restricted to Lie groups, see (1.4)

$$ G ~\subseteq~ {\cal D}iff^\infty ( M, M )$$

\medskip
and instead we can now deal with the {\it vastly larger} class of {\it semigroups}

$$ G ~\subseteq~ {\cal C}^\infty ( M, M ) $$

\medskip
In other words, in terms of the Euclidean domains $M$, we are expanding into the vastly larger
class of semigroups of smooth coordinate transforms which {\it need not} have smooth inverses,
and in fact, {\it need not} be surjective either. \\

Consequently, we are interested in exploring the structure of semigroups $G$ in
${\cal C}^\infty ( M, M )$ which have minimal overlap with ${\cal D}iff^\infty ( M, M )$. \\

Here we note two facts :

\begin{itemize}

\item Since we are interested in semigroups $G$ with neutral element $e \in G$, it follows
that a certain overlap between such semigroups and ${\cal D}iff^\infty ( M, M )$ is
inevitable. Indeed, in the overall semigroup ${\cal C}^\infty ( M, M )$, the neutral element
is $e = id_M$, that is, the identity mapping of $M$ into itself. And obviously, we have $e =
id_M \in {\cal D}iff^\infty ( M, M )$.

\item Similar with the classical Lie group theory, we may start with exploring the structure
of {\it one-dimensional} semigroups $G$ in ${\cal C}^\infty ( M, M )$.

\end{itemize}

In this regard, a first clarification follows from \\

{\bf Lemma 2.1} \\

Let $( X , \circ )$ be any semigroup with neutral element $e \in X$. Let

$$ [ 0, \infty ) \ni t \longmapsto x_t \in X $$

be a semigroup homomorphism, where $[ 0, \infty)$ is considered with its usual additive
semigroup structure. \\

If for a certain $t > 0$, the element $x_t \in X$ has an inverse in $X$, then every element
$x_s \in X$, with $s \in [ 0, \infty )$, has an inverse in $X$. \\
Consequently, the above semigroup homomorphism can be extended to a group homomorphism

$$ \mathbb{R} \ni t \longmapsto x_t \in X $$

\medskip
where $\mathbb{R}$ is considered with its usual additive group structure. \\

In particular, every element $x_t \in X$, with $t \in \mathbb{R}$, will have an inverse. \\

{\bf Proof} \\

Let $x^{\,\prime}_t \in X$ be the inverse element of $x_t$. Let $0 \leq s < t$, then the
commutativity of the additive semigroup on $[ 0, \infty )$ and the semigroup homomorphism
$[ 0, \infty ) \ni t \longmapsto x_t \in X$ sends $s + ( t - s)$ to $x_s \circ x_{t-s}$, and
$( t - s ) + s$ to $x_{t-s} \circ x_s$, both being equal with $x_t$, since $s + ( t - s) =
( t - s ) + s = t$. In this way

$$ x_s \circ x_{t-s} = x_t,~~~ x_{t-s} \circ x_s = x_t $$

Hence, by multiplying with $x^{\,\prime}_t $ on the right the first relation, and on the left
the second one, we obtain

$$ e = x_t \circ x^{\,\prime}_t = x_s \circ ( x_{t-s} \circ x^{\,\prime}_t ),~~~
   e = x^{\,\prime}_t \circ x_t = ( x^{\,\prime}_t \circ x_{t-s} ) \circ x_s $$

\medskip
Let us denote

$$ x^{\,\prime}_s = x_{t-s} \circ x^{\,\prime}_t,~~~
            x^{\,\prime\,\prime}_s = x^{\,\prime}_t \circ x_{t-s} $$

then

$$ x_s \circ x^{\,\prime}_s = x^{\,\prime\,\prime}_s \circ x_s = e $$

hence

$$ x^{\,\prime\,\prime}_s = x^{\,\prime\,\prime}_s \circ e = x^{\,\prime\,\prime}_s \circ
( x_s \circ x^{\,\prime}_s ) = ( x^{\,\prime\,\prime}_s \circ x_s ) \circ x^{\,\prime}_s =
e \circ x^{\,\prime}_s = x^{\,\prime}_s$$

\medskip
Therefore $x_s \in X$ has indeed the inverse $x^{\,\prime}_s = x^{\,\prime\,\prime}_s
\in X$. \\

Now, let $t < s < \infty$. Let $n \in \mathbb{N}$ be such that $s < n t$. Then obviously
$x_{n t} \in X$ has the inverse $x^{\,\prime}_t \circ ~.~.~.~ \circ x^{\,\prime}_t \in X$,
where the composition $\circ$ is applied $n-1$ times. Taking now $t^{\,\prime} = n t$, we
obtain $s < t^{\,\prime}$, and the proof is reduced to the previous case. \\

{\bf Remark 2.1} \\

In the above Lemma it is not necessary that the semigroup $( X, \circ )$ be commutative, nor
that $x_0 = e \in X$.

\hfill $\Box$ \\

In view of Lemma 2.1, we are led to, Rosinger [1, chap. 13] \\

{\bf Definition 2.1 ( Genuine Lie Semigroups, or GLS )} \\

We call {\it one-dimensional genuine Lie semigroup} on $M$, every semigroup homomorphism \\

(2.1) $~~~ [ 0, \infty ) \ni t \longmapsto g_t \in {\cal C}^\infty ( M, M ) $ \\

which has the properties \\

(2.2) $~~~ g_0 ~=~ id_M $ \\

(2.3) $~~~ g_t \in {\cal C}^\infty ( M, M ) \setminus {\cal D}iff^\infty ( M, M ),
                                                               ~~~\mbox{for}~~ t > 0 $ \\

{\bf Remark 2.2} \\

1) Actually, every genuine Lie semigroup, or GLS, in the above definition is given by the
image of the semigroup homomorphism (2.1), namely \\

(2.4) $~~~ G ~=~ \{~ g_t ~~|~~ t \in [ 0, \infty ) ~\} ~\subseteq~
                                               {\cal C}^\infty ( M, M ) $ \\

And we call it {\it genuine}, since we obviously have \\

(2.5) $~~~ G ~\cap~ {\cal D}iff^\infty ( M, M ) ~=~ \{~ id_M ~\} $ \\

In other words, {\it none} of the elements of $G$, except for $id_M$, are invertible in the
overall semigroup ${\cal C}^\infty ( M, M )$. \\

2) In view of Lemma 2.1, we have the following {\it dichotomy} with respect to arbitrary
semigroup homomorphisms \\

(2.6) $~~~ [ 0, \infty ) \ni t \longmapsto g_t \in  {\cal C}^\infty ( M, M ),~~~
                                                        \mbox{with}~~ g_0 = id_M $ \\

namely : \\

Either

\begin{itemize}

\item None of the smooth coordinate transforms $g_t$ has a smooth inverse, except for $g_0 =
id_M$.

\end{itemize}

Or

\begin{itemize}

\item All the smooth coordinate transforms $g_t$ have a smooth inverse, and then the above
semigroup homomorphism (2.6) can be extended to a group homomorphism

\end{itemize}

(2.7) $~~~ \mathbb{R} \ni t \longmapsto g_t \in  {\cal D}iff^\infty ( M, M ) $ \\

In view of this dichotomy, and within the realm of semigroups of transformations, the concept
of {\it genuine} one-dimensional Lie semigroup proves to be the {\it natural} alternative to
that of one-dimensional Lie group. \\

3) Here it should be recalled that Lie semigroups have been studied in Hilgert et.al., for
instance. However, so far, such studies have only concerned Lie semigroups which are
sub-semigroups of Lie groups, or with the above notation, are sub-semigroups of ${\cal
D}iff^\infty ( M, M )$. \\
Therefore, they {\it cannot} be genuine Lie semigroups.

\hfill $\Box$ \\

In the sequel, it will be shown that there are {\it plenty} of one-dimensional genuine Lie
semigroups. Namely, it will among others be shown that \\

(2.8) $~~~ \begin{array}{l}
              \forall~~ f \in {\cal C}^\infty ( M, M ) : \\ \\
              \exists~~ [ 0, \infty ) \ni t \longmapsto g_t \in {\cal C}^\infty ( M, M )
                       ~~\mbox{one-dimensional GLS} : \\ \\
               ~~~~ f ~=~ g_1
            \end{array} $ \\

In other words, when taken all together, the one-dimensional genuine Lie semigroups {\it
cover} the whole of ${\cal C}^\infty ( M, M )$. \\
Needless to say, any given one-dimensional genuine Lie semigroup is but a {\it path} in
${\cal C}^\infty ( M, M )$. \\ \\

{\bf 3. Applications to PDEs, Semi-symmetries} \\

Very simple examples can already show that {\it genuine Lie semigroups}, or in short, GLS, can
give {\it semi-symmetries} of PDEs which cannot be obtained by smooth coordinate transforms
which are invertible. Therefore, they are {\it not} within the reach of Lie group theory. \\

Let us start with an example, before giving the precise definition. \\

One of the simplest linear PDEs is \\

(3.1) $~~~ U_t = U_x,~~~ ( t, x ) \in \Omega = \mathbb{R}^2 $ \\

whose smooth solutions are given by \\

(3.2) $~~~ U ( t, x ) = h ( t + x ),~~~ ( t, x ) \in \Omega $ \\

where $h \in {\cal C}^\infty ( \mathbb{R}, \mathbb{R} )$ is arbitrary. \\

Let us take now $M = \Omega \times \mathbb{R} = \mathbb{R}^3$, and any smooth function $f \in
{\cal C}^\infty ( M, M )$ which is of the form \\

(3.3) $~~~ M \ni ( t, x, u ) \stackrel{f}\longmapsto ( t, x, g ( u ) ) \in M $ \\

where $g \in {\cal C}^\infty ( \mathbb{R}, \mathbb{R} )$ is an arbitrary {\it non-injective}
function. \\
Then clearly, the smooth function $f : M \longrightarrow M$ is {\it not} a coordinate
transform on $M$, since it is not injective, thus it {\it fails} to be invertible, and
consequently $f \notin {\cal D}iff^\infty ( M, M)$. \\

We note, nevertheless, that the function $f$ {\it transforms} solutions of the PDE in (3.1)
into solutions of the same PDE, thus it is a {\it semi-symmetry} of that equation. Indeed, the
action of the function $f$ upon the solution $U$ in (3.2) is given by,  Appendix \\

(3.4) $~~~ \widetilde{U} ( t, x ) = ( f U ) ( t, x ) = g ( h ( t + x ) ),~~~
                                                        ( t, x ) \in \Omega $ \\

thus $\widetilde{U} = f U$ is again a solution of the PDE in (3.1). \\

The above example of a semi-symmetry which is {\it not} a usual symmetry is no doubt very
simple. However, it can already clearly illustrate the main issue. Namely, given any linear or
nonlinear PDE \\

(3.5) $~~~ T ( x, D )~ U ( x ) ~=~ 0,~~~ x \in \Omega $ \\

where $\Omega$ is an open subset in $\mathbb{R}^n$, and $U : \Omega \longrightarrow
\mathbb{R}$ is the unknown solution. If we denote $M = \Omega \times \mathbb{R}$, then the
standard Lie symmetries of that equation are given by functions $f \in {\cal D}iff^\infty
( M, M )$ which by their actions, turn solutions $U$ into solutions $\widetilde{U} = f U$ of
that equation. \\

However, as seen above in (3.4), we can define such actions not only for functions $f \in
{\cal D}iff^\infty ( M, M )$, but also for the {\it much larger} class of functions $f \in
{\cal C}^\infty ( M, M )$. And such actions $\widetilde{U} = f U$ are called {\it
semi-symmetries} of the PDEs in (3.5), if they turn solutions $U$ of those equations into
solutions $\widetilde{U} = f U$ of the respective equations. \\

Here it should be noted that in order to define the actions $\widetilde{U} = f U$ for
arbitrary solutions $U \in {\cal C}^\infty ( M, \mathbb{R} )$ and functions $f \in {\cal
C}^\infty ( M, M )$, one has to use the full power of the {\it parametric} representation of
functions and actions, introduced and developed in Rosinger [1, chapters 1-5]. For
convenience, a brief review of these issues is presented in the Appendix. \\

Also, in Rosinger [1, chap. 13], other more involved examples of semi-symmetries of PDEs are
presented. \\ \\

{\bf 4. "Enforcing" : How to Generate One-Dimensional Genuine \\
\hspace*{0.7cm} Lie
Semigroups} \\

Our interest is to generate one-dimensional genuine Lie semigroups in ${\cal C}^\infty ( M,
M )$, for arbitrary nonvoid open sets $M \subseteq \mathbb{R}^l$. \\

What we shall do in this section is to show how to solve the following particular case of that
general problem. Namely, let us take arbitrary smooth functions \\

$~~~~~~ f \in {\cal C}^\infty ( M, M ) $ \\

and use them in a simple "enforcing" method, in order to generate one-dimensional genuine Lie
sub-semigroups \\

$~~~~~~ G^\# ~\subset~ {\cal C}^\infty ( M^\#, M^\# ) $ \\

such that \\

$~~~~~~ f^\# ~=~ ( 1, f ) \in G^\# $ \\

where \\

$~~~~~~ M^\# ~=~ ( 0, \infty ) \times M ~\subseteq~ \mathbb{R}^{l + 1} $ \\

In this way, we obtain the following general result \\

( dim + 1 ) $~~~~~~ \begin{array}{l}
                      \forall~~ M \subseteq \mathbb{R}^l ~~\mbox{nonvoid open},~
                                    f \in {\cal C}^\infty ( M, M ) ~: \\ \\
                      \exists~~ G^\# ~\subset~ {\cal C}^\infty
                                  ( M^\#, M^\# ) ~~\mbox{a genuine Lie semigroup} ~: \\ \\
                      ~~~~  f^\# ~=~ ( 1, f ) \in G^\#
                   \end{array} $ \\ \\

where $M^\# ~=~ ( 0, \infty ) \times M ~\subseteq~ \mathbb{R}^{l + 1}$. \\

{\bf Open Problem.} The ultimately general result regarding GLS-s, namely \\

( GLS ) $~~~~~~ \begin{array}{l}
                      \forall~~ M \subseteq \mathbb{R}^l ~~\mbox{nonvoid open},~
                                    f \in {\cal C}^\infty ( M, M ) ~: \\ \\
                      \exists~~ G ~\subset~ {\cal C}^\infty ( M, M )
                                  ~~\mbox{a genuine Lie semigroup} ~: \\ \\
                      ~~~~  f \in G
                   \end{array} $ \\

is still open. \\

{\bf Short Review of the Classical Case of Lie Groups.} When generating one-dimensional
genuine Lie semigroups, we shall try as much as possible to follow the classical way in Lie
group theory which generates the one-dimensional Lie groups. There are a number of well known
reasons why that classical way is important and useful, and therefore, its possible extension
to the generation of one-dimensional genuine Lie semigroups is worth exploring. Indeed, in
that classical way, several fundamental mathematical ideas, constructions and properties come
together in a fruitful interaction. Among them are :

\begin{itemize}

\item Lie Groups,

\item Actions,

\item Infinitesimal Generators,

\item Autonomous ODEs,

\item Flows,

\item Evolution Operators.

\end{itemize}

Let us start, therefore, by recalling here in short that classical
construction. \\

{\bf Autonomous ODEs Generating Lie Groups.} Let $M \subseteq \mathbb{R}^l$ be any nonvoid
open subset and $F \in {\cal C}^\infty ( M, \mathbb{R}^l )$ any smooth function. We consider
the {\it autonomous} nonlinear system of ODEs with the respective initial conditions \\

(4.1) $~~~ \begin{array}{l}
                   \frac{\D d}{\D dt} Y ( t ) = F ( Y ( t ) ),~~~ t \in \mathbb{R} \\ \\
                   Y ( 0 ) = y \in M
           \end{array} $ \\

Then the unique solution - for convenience, assumed to exist globally - namely \\

(4.2) $~~~ \mathbb{R} \ni t \longmapsto Y ( t ) \in M $ \\

defines through each point $y \in M$ a {\it flow} on $M$, which corresponds to the {\it
one-dimensional Lie group}~ $G = ( \mathbb{R}, + )$ that {\it acts}~ on $M$ according to \\

(4.3) $~~~ \mathbb{R} \times M \ni ( t, y ) \longmapsto Y ( t ) \in M $ \\

In this case $F$ is called the {\it infinitesimal generator} of the Lie group action (4.3). \\

{\bf Lie Groups Generating Autonomous ODEs.} A basic fact in Lie group theory is that the
converse of the above construction also holds. Namely, given the one-dimensional Lie group
action (4.3) on $M$, then one can obtain an infinitesimal generator $F \in {\cal C}^\infty
( M, \mathbb{R}^l )$ defined by \\

(4.4) $~~~ F ( y ) = \frac{\D d}{\D dt} Y ( t ) |_{t = 0},~~~ y \in M $ \\

and the corresponding ODE and initial value problem (4.1) will always have a global solution.
In this case, the steps (4.1) - (4.3) will give us back the initial Lie group action (4.3) on
$M$, with which we started. \\

{\bf Evolution Operators.} The above in (4.1) - (4.4) can be described in terms of {\it
evolution operators} $E$ as well. Namely, we define \\

(4.5) $~~~ \begin{array}{l}
                \mathbb{R} \ni t \longmapsto E ( t ) : M \longrightarrow M \\ \\
                E ( t ) ( y ) = Y ( t ),~~~ t \in \mathbb{R},~ y \in M
           \end{array} $ \\

and then the above one-dimensional Lie group action (4.3) can be written in the form \\

(4.6) $~~~ \mathbb{R} \times M \ni ( t, y ) \longmapsto E ( t ) ( y ) = Y ( t ) \in M $ \\

These evolution operators $E$ have the important {\it group} property \\

(4.7) $~~~ \begin{array}{l}
               E ( 0 ) = id_M \\ \\
               E ( t ) \circ E ( s ) = E ( t + s ),~~~ t, s \in \mathbb{R}
            \end{array} $ \\

thus we have the {\it group homomorphism}, see (4.5) \\

(4.8) $~~~ \mathbb{R} \ni t \longmapsto E ( t ) \in {\cal D}iff^\infty ( M, M) $ \\

{\bf "Enforcing" : How to Find Infinitesimal Generators for One-Dimensional Genuine Lie
Semigroups.} Our aim - according to the most general program ( GLS ) above - is to find on $M$
one-dimensional semigroup actions, see (4.3), (4.6) \\

(4.9) $~~~ [ 0, \infty ) \times M \ni ( t, y ) \longmapsto S ( t, y ) \in M $ \\

or equivalently \\

(4.10) $~~~ [ 0, \infty ) \ni t \longmapsto
                    S ( t, . ) : M \ni y \longmapsto S ( t, y ) \in M $ \\

which give one-dimensional genuine Lie semigroups, that is, with the properties, see (2.1) -
(2.3) \\

(4.11) $~~~ \begin{array}{l}
                S ( 0, . ) = id_M \\ \\
                S ( t, . ) \circ S ( s, . ) = S ( t + s, . ),~~~ t,~ s \in [ 0, \infty ) \\ \\
                S ( t, . ) \notin {\cal D}iff^\infty ( M, M ),~~~ t > 0
            \end{array} $ \\

For that purpose, and in order to become more familiar with the new one-dimensional genuine
Lie semigroup situation, we shall proceed step by step, analyzing cases which are more and
more general, and in the process we shall eliminate those which are not suited. \\

{\bf Autonomous Nonsingular ODEs with Global Solutions.} First we note that a one-dimensional
GLS in (4.9) - (4.11) {\it cannot} be generated by an autonomous ODE of the type (4.1). \\
Indeed, if the respective ODE in (4.1) has global solutions (4.2) for every $y \in M$, then as
seen above, it generates a one-dimensional Lie group (4.8) acting on $M$, which obviously is
{\it not} a genuine Lie semigroup. \\

{\bf Autonomous Nonsingular ODEs with Local Solutions.} On the other hand, in case the ODE in
(4.1) does not have such a global solution property, then since the respective $F$ is assumed
to be smooth on the whole of $M$, a classical result on the existence of solutions for ODEs,
Coddington \& Levinson, states that for every initial condition $y \in M$, there exists a
largest nonvoid open interval $0 \in I_y \subseteq \mathbb{R}$, such that a unique solution $Y
: I_y \longrightarrow M$ exists, which satisfies the initial condition $Y ( 0 ) =y \in M$. \\
Here it is important to note that such a unique solution $Y$ will exist on an open
neighbourhood of $t = 0 \in \mathbb{R}$, that is, both for strictly positive and strictly
negative values of $t$, possibly limited accordingly. And as seen next, this is enough in
order to prevent such a solution $Y$ from generating a genuine Lie semigroup on $M$. \\
Indeed, the above unique solution property gives \\

(4.12) $~~~ M \ni y \longmapsto E : I_y \ni t \longmapsto E ( t ) ( y ) = Y ( t ) \in M $ \\

Let us denote for $t \in \mathbb{R}$ \\

(4.13) $~~~ M_t ~=~ \{~ y \in M ~~|~~ t \in I_y ~\} $ \\

then (4.12) gives \\

(4.14) $~~~ \mathbb{R} \ni t \longmapsto E ( t ) : M_t \ni y \longmapsto
                                                E ( t ) ( y ) = Y ( t ) \in M $ \\

and we have the {\it generalized} group property of the evolution operators $E$, Rosinger [3,
pp. 56,57], given by the commutative diagram \\

\begin{math}
\setlength{\unitlength}{0.2cm}
\thicklines
\begin{picture}(60,20)

\put(10,16){$M_{t+s}$}
\put(27,18){$E ( t + s )$}
\put(15,16.5){\vector(1,0){31}}
\put(48,16){$M_0 = M$}
\put(0,6){$(4.15)$}
\put(15,14.5){\vector(1,-1){13.5}}
\put(16,6){$E ( t )$}
\put(30,0){$M_s$}
\put(34,1){\vector(1,1){13.5}}
\put(42,6){$E ( s )$}

\end{picture}
\end{math} \\

where $t, s \in \mathbb{R}$. \\
What happens now in case the ODE in (4.1) does not have global solutions (4.2) for every $y
\in M$, is that $M_t = \phi$, or at least $M_t \subset M,~ M_t \neq M$, for certain $t \in
\mathbb{R}$. Consequently, the evolution operators $E$ are {\it not} defined on the whole of
$M$, thus we cannot possibly obtain (4.10), where $S ( t, . )$ are supposed to be defined
everywhere on $M$. \\

{\bf Autonomous Singular ODEs.} A next level of generality is to consider the autonomous ODEs
in (4.1) with $F$ no longer smooth all over $M$, but having certain {\it singularities}, for
instance \\

(4.16) $~~~ F \in {\cal C}^\infty ( M \setminus \Sigma,~ \mathbb{R}^l )
                                           \setminus {\cal C}^\infty ( M, \mathbb{R}^l ) $ \\

for suitable nonovid subsets $\Sigma \subset M$. Such an approach, however, need {\it not}
always lead to one-dimensional genuine Lie semigroups, as seen from the following simple
example. Let $M = \mathbb{R}$, and consider the ODE \\

$~~~~~~ \begin{array}{l}
              \frac{\D d}{\D dt} Y ( t ) = 1 / ( Y ( t ) )^2,~~~ t \in \mathbb{R} \\ \\
              Y ( 0 ) = y \in M \setminus \{ 0 \}
        \end{array} $ \\

Here we have $F ( y ) = 1 / y^2$, for $y \in M \setminus \Sigma$, where $\Sigma = \{ 0 \}
\subset M$, thus (4.16) is satisfied. However, the unique solution is \\

$~~~~~~ Y ( t ) = ( 3 t + y^3 )^{1/3},~~~ y \in M,~ y \neq 0,~
                               t \in \mathbb{R},~ t \neq - y^3 / 3 $ \\

And we note that this function $Y$ can in fact be extended to \\

$~~~~~~ Y ( t ) = ( 3 t + y^3 )^{1/3},~~~ y \in M,~ t \in \mathbb{R} $ \\

which for every $y \in M$ satisfies the singular autonomous ODE \\

$~~~~~~ \frac{\D d}{\D dt} Y ( t ) = 1 / ( Y ( t ) )^2,~~~
                                  t \in \mathbb{R},~ t \neq - y^3 / 3 $ \\

Thus instead of (4.10), we have \\

$~~~~~~ \mathbb{R} \ni t \longmapsto S ( t, . ) : M \ni y \longmapsto Y ( t ) \in M $ \\

and this leads to a Lie group action on $M$, since obviously \\

$~~~~~~ S ( t, S ( s, y ) ) = S ( t + s, y ),~~~ t, s \in \mathbb{R},~ y \in M $ \\

{\bf Non-autonomous Singular ODEs.} In view of the above, the next step is to consider {\it
non-autonomous} ODEs of the form \\

(4.17) $~~~ \begin{array}{l}
               \frac{\D d}{\D dt} Y ( t ) = F ( t, Y ( t ) ),~~~ t \in \mathbb{R} \\ \\
               Y ( t_0 ) = y_0 \in M
            \end{array} $ \\

where $F : \mathbb{R} \times M \longrightarrow \mathbb{R}^l$ and \\

(4.18) $~~~ F \in {\cal C}^\infty $ \\

except for certain possible singularities in its domain $\mathbb{R} \times M$. \\

The idea here is twofold, namely :

\begin{itemize}

\item {\bf ( Reduct )}~~ To use the standard {\it reduction} method of such non-autonomous
ODEs to autonomous ones.

\item {\bf ( Sing )}~~ To include certain {\it singularities}\, in the non-autonomous ODEs, so
that, when reduced to autonomous ODEs, the solutions of these latter ODEs do {\it not} give
one-dimensional Lie groups, but only one-dimensional genuine Lie semigroups.

\end{itemize}

It follows that the only problem here is to find out what kind of {\it singularities}\, the
non-autonomous ODEs (4.17), and more precisely, their right hand terms (4.18), must have in
the very least, in order to secure the above property {\bf ( Sing )}. \\

{\bf Remark 4.1.} \\

For the sake of clarity, let us recall in short the standard way non-autonomous ODEs can be
{\it reduced} to autonomous ones. The further details needed will be presented in section
5. \\
Given an explicit non-autonomous ODE with a respective initial value problem \\

(4.19) $~~~ \begin{array}{l}
                   \frac{\D d}{\D dt} Y ( t ) = F ( t, Y ( t ) ),~~~ t \in \mathbb{R} \\ \\
                   Y ( t_0 ) = y_0 \in M
           \end{array} $ \\

or more generally, an implicit non-autonomous ODE with an associated initial value problem \\

(4.20) $~~~ \begin{array}{l}
                   F ( t, Y ( t ), \frac{\D d}{\D dt} Y ( t ) ) = 0,~~~ t \in \mathbb{R} \\ \\
                   Y ( t_0 ) = y_0 \in M
           \end{array} $ \\

there is a well known standard procedure in Control Theory to {\it reduce} it to an autonomous
ODE. This is done simply by {\it increasing} with 1 the dimension of the system of ODEs (4.20),
namely, from $l$ to $l + 1$. For that purpose, we replace the $l$-dimensional {\it solution}
vector $Y : \mathbb{R} \longrightarrow M$ with the $(l + 1)$-dimensional solution vector \\

(4.21) $~~~ Y^\# : \mathbb{R} \longrightarrow M^\# $ \\

where \\

(4.22) $~~~ M^\# = \mathbb{R} \times M,~~~
            Y^\# ( t ) = ( t, Y ( t ) ),~~~ t \in \mathbb{R} $ \\

In this case (4.20) obviously becomes the implicit autonomous ODE \\

(4.23) $~~~ \begin{array}{l}
                   F^\# ( Y^\# ( t ), \frac{\D d}{\D dt} Y^\# ( t ) ) = 0,~~~
                                                t \in \mathbb{R} \\ \\
                   Y^\# ( t_0 ) = ( t_0, y_0 ) \in M^\#
           \end{array} $ \\

where the equation $F^\# ( Y^\# ( t ), \frac{\D d}{\D dt} Y^\# ( t ) ) = 0$ given by \\

(4.24) $~~~ \begin{array}{l}
                  \frac{\D d}{\D dt}\, t = 1,~~~ t \in \mathbb{R} \\ \\
                  F ( t, Y ( t ), \frac{\D d}{\D dt} Y ( t ) ) = 0,~~~ t \in \mathbb{R} \\ \\
            \end{array} $ \\

Clearly, in the particular case of (4.19), this procedure leads to an autonomous ODE system
which again is explicit, namely \\

(4.23$^*$) $~~~ \frac{\D d}{\D dt} Y^\# ( t ) = F^\# ( Y^\# ( t ) ),~~~ t \in \mathbb{R} $ \\

However, it is important to note that in (4.23) care has to be taken with the {\it initial
condition}. Indeed, if as in (4.20), the initial condition is given at $t_0 \in \mathbb{R}$,
then in the case of the extended ODE in (4.23), this will become \\

(4.25) $~~~  Y^\# ( t_0 ) = ( t_0, y_0 ) \in M^\# $ \\

in other words, the right hand term in (4.25) is {\it not} completely arbitrary in $M^\#$,
since it {\it must} have the {\it same}\, $t_0$ as in the left hand term. \\

Now the point of interest for us in the above reduction of non-autonomous ODEs to autonomous
ones is in the following two facts :

\begin{itemize}

\item The resulting autonomous ODE in (4.23) has solutions $Y^\#$ defined on the {\it same}
$t$-interval $I \subseteq \mathbb{R}$ with the solutions $Y$ of the non-autonomous ODE (4.20),
as this follows easily from (4.24).

\item The solutions $Y^\#$ of the autonomous ODE in (4.23) may under rather general conditions
have a one-dimensional {\it group} property similar with (4.5) - (4.7).

\item In case there are appropriate {\it singularities}\, involved in the ODE in (4.23), the
solutions $Y^\#$ will have a one-dimensional {\it semigroup} property, namely, they lead to an
evolution operator

\end{itemize}

(4.26) $~~~ \begin{array}{l}
            [ 0, \infty ) \ni t \longmapsto E^\# ( t ) : M^\# \longrightarrow M^\# \\ \\
            E^\# ( t ) ( y^\# ) = Y^\# ( t ),~~~ t \in [ 0, \infty ),~ y^\# \in M^\#
           \end{array} $ \\

\begin{quote}

with the semigroup property \\

\end{quote}

(4.27) $~~~ \begin{array}{l}
               E^\# ( 0 ) = id_{M^\#} \\ \\
               E^\# ( t ) \circ E^\# ( s ) = E^\# ( t + s ),~~~ t, s \in [ 0, \infty )
            \end{array} $ \\

In this way, we shall be able to generate {\it one-dimensional genuine Lie semigroups} in
${\cal C}^\infty ( M^\#,~ M^\# )$, rather than in the initial ${\cal C}^\infty
( M,~ M )$. \\

{\bf A First Simple Example of "Enforcing" : A Singular Non-autonomous ODE and its Solution.}
We note that we can start, so to say, {\it backwards}. Namely, we can start with an {\it
action} on $M$ of the form, see (4.9) - (4.10), given by \\

(4.28) $~~~ [ 0, \infty ) \times M \ni ( t, y ) \longmapsto H ( t, y ) \in M $ \\

which need not be a semigroup, and then find the non-autonomous singular ODE which it
satisfies. \\

Now the {\it enforcing} part comes when on purpose we choose this action (4.28) so that it
{\it cannot} be extended to an action \\

(4.28$^*$) $~~~ ( - \infty, \infty ) \times M \ni ( t, y ) \longmapsto H ( t, y ) \in M $ \\

which in certain cases may possibly give a Lie group action on $M$. \\

Consequently, by such an enforcing, we may obtain an insight into the kind of {\it
singularities} possessed by the non-autonomous ODE satisfied by such an action (4.28). \\

As mentioned, such an action (4.28) from which we start need not in general be a semigroup
action on $M$, and even less a genuine Lie semigroup action. However, when the respective
non-autonomous ODE satisfied by that action is reduced to an autonomous one, we shall
inevitably obtain a semigroup action, this time on the $(l+1)$-dimensional open subset $M^\#
= \mathbb{R} \times M$, see Remark 4.1, above. \\

And now, we choose one of the simplest possible such actions (4.28), namely \\

(4.29) $~~~ H ( t, y ) = y + \sqrt t\, y^2,~~~ t \in [ 0, \infty ),~ y \in M = \mathbb{R} $ \\

Then obviously (4.29) cannot be extended to a smooth action (4.28$^*$) owing to the presence
of $\sqrt t$. Also, since $y^2$ appears in (4.29), the corresponding action (4.28) is not
injective. Furthermore, we have \\

(4.30) $~~~ \begin{array}{l}
                   H ( 0, y ) = y,~~~ y \in M \\ \\
                   H ( t, . ) \notin {\cal D}iff^\infty ( M, M ),~~~ t \in ( 0, \infty ) \\ \\
                   H \notin {\cal C}^1 ( [ 0, \infty ) \times M, M ) \\ \\
                   H \in {\cal C}^0 ( [ 0, \infty ) \times M, M ) \cap
                         {\cal C}^\infty ( ( 0, \infty ) \times M, M )
             \end{array} $ \\

Let us now find the non-autonomous explicit or implicit ODE satisfied by $H$ in (4.29). By
partial derivation of that relation with respect to $t$, we obtain \\

(4.31) $~~~ \frac{\D \partial}{\D \partial t} H ( t, y ) = \frac{\D y^2}{\D 2 \sqrt t},~~~
                                                        t \in ( 0, \infty ),~ y \in M $ \\

Then we solve (4.29) as a quadratic equation in $y$, and obtain \\

(4.32) $~~~ \begin{array}{l}
            y_1,~ y_2 =
                    \frac{\D -1 \pm \sqrt{1 + 4 \sqrt t\, H ( t, y )}}{\D 2 \sqrt t}~, \\
~~~~~~~~~~~~~~~~~~~~~~~~~~~~~~~~~~~ t \in ( 0, \infty ),~ 1 + 4 \sqrt t\, H ( t, y ) \geq 0
            \end{array} $ \\

Further, we return to (4.29), replace $y$ with its values from (4.32), while we replace $y^2$
with its value from (4.31). In this way, we obtain the two {\it non-autonomous singular
explicit} ODEs \\

(4.33) $~~~ \begin{array}{l}
          \frac{\D \partial}{\D \partial t} H ( t, y ) =
            \frac{\D 1 + 2 \sqrt t\, H ( t, y ) \pm \sqrt{1 + 4 \sqrt t\, H ( t, y )}}
                                               {\D 4 t \sqrt t}~,~~~ \\ \\
~~~~~~~~~~~~~~~~~~~~~~~~~~~~~~~~~~~ t \in ( 0, \infty ),~ 1 + 4 \sqrt t\, H ( t, y)  \geq 0
             \end{array} $ \\

If we now consider $H$ given in (4.29) as a function of $t \in ( 0, \infty) $, and with $y \in
M$ being a fixed parameter, then substituting $H$ into (4.33), a simple computation shows that
$H$ satisfies the ODE in (4.33) with the sign "$+$" in front of the large radical, when \\

(4.34) $~~~ t \in ( 0, \infty ),~ y \in M,~ 1 + 2 \sqrt t\, y \leq 0 $ \\

and alternatively, satisfies the ODE in (4.33) with the sign "$-$" in front of the large
radical, when \\

(4.35) $~~~ t \in ( 0, \infty ),~ y \in M,~ 1 + 2 \sqrt t\, y \geq 0 $ \\

Moreover, in both of these cases, $H$ as defined in (4.29) will obviously satisfy the {\it
limit} type initial condition \\

(4.36) $~~~ \lim_{_{~t~ \searrow~ 0}}~ Y ( t ) = y $ \\

since in fact we have $H ( 0, y ) = y$, for $y \in M$, as for every $y \in M$, H is continuous
at $t = 0$ from the right, see (4.30). \\

Here we note that, since it is necessary that $t > 0$ for the ODEs in (4.33) to be defined,
one cannot in general ask for them to satisfy a usual initial condition $Y ( 0 ) = y$, but
only of the limit type one in (4.36). \\

{\bf A Type of Implicit Singular Non-autonomous ODEs Leading to Genuine Lie Semigroups.} In
view of the above example, we are led to consider {\it explicit} non-autonomous singular ODEs
and associated limit type initial condition of the form, see (4.17), (4.18) \\

(4.37) $~~~ \begin{array}{l}
               \frac{\D d}{\D dt} Y ( t ) = F ( t, Y ( t ) ),~~~ t \in ( 0, \infty ) \\ \\
               \lim_{_{~t~ \searrow~ 0}}~ Y ( t ) = y \in M_0
            \end{array} $ \\

for a suitable subset $M_0 \subseteq M$, and with \\

(4.38) $~~~ F \in {\cal C}^\infty ( ( 0, \infty ) \times M,~ \mathbb{R}^l ) $ \\

And we expect to have unique solutions $Y$ such that \\

(4.39) $~~~ \begin{array}{l}
            ~~~*)~~ Y \in {\cal C}^\infty ( ( 0, \infty ), M ) \\ \\
            ~**)~~ Y \in {\cal C}^0 ( [ 0, \infty ), M )
                                     \setminus {\cal C}^1 ( [ 0, \infty ), M ) \\ \\
            ~***)~~ ( 0, \infty ) ~~~\mbox{is the largest interval in}~~ \mathbb{R} \\
                     ~~~~~~~~~~~\mbox{with properties}~~ *) ~~\mbox{and}~~ **)
            \end{array} $ \\

Needless to say, we may have to go one level of generality higher, namely, to consider {\it
implicit} non-autonomous singular ODEs and associated limit type initial conditions, of the
form \\

(4.40) $~~~ \begin{array}{l}
                   F ( t, Y ( t ), \frac{\D d}{\D dt} Y ( t ) ) = 0,~~~
                                                t \in ( 0, \infty ) \\ \\
                   \lim_{_{~t~ \searrow~ 0}}~ Y ( t ) = y \in M_0
            \end{array} $ \\

instead of the explicit ones in (4.37). \\

We present a first result on the existence of solutions (4.39) of implicit singular
non-autonomous ODEs of the type (4.40). This result shows that there are so many such
solutions as to {\it cover} the whole of ${\cal C}^\infty ( M, M )$, at least in the case when
$M = \mathbb{R}^l$. \\

{\bf Lemma 4.1. "Enforcing"} \\

Let $f \in {\cal C}^\infty ( \mathbb{R}^l, \mathbb{R}^l ) \setminus {\cal D}iff^\infty
( \mathbb{R}^l, \mathbb{R}^l )$. Then there exist actions on $\mathbb{R}^l$ \\

(4.41) $~~~ \begin{array}{l}
                [ 0, \infty ) \times \mathbb{R}^l \ni ( t, y ) \longmapsto
                                                       H ( t, y ) \in \mathbb{R}^l \\ \\
                H \in {\cal C}^0 ( [ 0, \infty ) \times \mathbb{R}^l,~ \mathbb{R}^l ) \cap
                   {\cal C}^\infty ( ( 0, \infty ) \times \mathbb{R}^l,~ \mathbb{R}^l )
             \end{array} $ \\

such that \\

(4.42) $~~~ H ( 0, . ) = id_{\mathbb{R}^{\,l}},~~~ H ( 1, . ) = f $ \\

and $H$ as a function of $t$, and with $y \in \mathbb{R}^l$ considered as a fixed parameter,
is a solution of an implicit singular non-autonomous ODE and associated limit type initial
condition in (4.40), which satisfies (4.39). \\

{\bf Proof} \\

Let us take any function $g \in {\cal C}^0 ( [ 0, \infty ), \mathbb{R} ) \cap {\cal C}^\infty
( ( 0, \infty ), \mathbb{R} )$ such that \\

(4.43) $~~~ \begin{array}{l}
                  g ( 0 ) = 0,~~~ g ( 1 ) = 1 \\ \\
                  g \notin {\cal C}^1 ( [ 0, \infty ), \mathbb{R} ) \\ \\
                  g^{\,\prime} ( t ) \neq 0,~~~ t \in ( 0, \infty )
             \end{array} $ \\

for instance, we can consider $g ( t ) = \sqrt t$, for $t \in [ 0, \infty )$. \\

Let us then define \\

(4.44) $~~~  H \in {\cal C}^0 ( [ 0, \infty ) \times \mathbb{R}^l,~ \mathbb{R}^l ) \cap
                   {\cal C}^\infty ( ( 0, \infty ) \times \mathbb{R}^l,~ \mathbb{R}^l )$ \\

as the homotopic deformation of $id_{\mathbb{R}^{\,l}}$ into $f$, mediated by $g$, namely \\

(4.45) $~~~ H ( t, y ) = ( 1 - g ( t ) ) y + g ( t ) f ( y ),~~~
                               t \in [ 0, \infty ),~ y \in \mathbb{R}^l $ \\

Then clearly \\

(4.46) $~~~ H ( 0, . ) = id_{\mathbb{R}^{\,l}},~~~ H ( 1, . ) = f \in {\cal C}^\infty
                               ( \mathbb{R}^l, \mathbb{R}^l ) \setminus
                                   {\cal D}iff^\infty ( \mathbb{R}^l, \mathbb{R}^l )$ \\

Let us now find the ODE satisfied by $H$ when it is considered a function of $t$, while $y$ is
taken as a parameter. The relation (4.45) and its partial derivative in $t$ give the two
linear algebraic equations in $y$ and $f ( y )$, namely \\

(4.47) $~~~ \begin{array}{l}
                  H = ( 1 - g ( t ) ) y + g ( t ) f ( y ) \\ \\
                  \frac{\D \partial}{\D \partial t} H = g^{\,\prime} ( t ) ( f ( y ) - y )
             \end{array} $ \\

for $t \in ( 0, \infty ),~ y \in M$. Since in view of (4.43) we have $g^{\,\prime} ( t ) \neq
0$ for $t \in ( 0, \infty )$, we can solve (4.47) for $y$ and $f ( y )$, and obtain \\

(4.48) $~~~ y = \frac{\D g^{\,\prime} H - g H_t }{\D g^{\,\prime} },~~~
            f ( y ) = \frac{\D ( 1 - g ) H_t + g^{\,\prime} h}{\D g^{\,\prime} },~~~
                             t \in ( 0, \infty ) $ \\

Substituting in (4.45) these values for $y$ and $f ( y )$, we obtain the {\it implicit
singular non-autonomous} ODE in $Y$, namely \\

(4.49) $~~~ ( 1 - g ( t ) ) \frac{\D d}{\D dt} Y ( t ) + g^{\,\prime} ( t ) Y ( t ) =
                g^{\,\prime} ( t ) f \left ( \frac{\D g^{\,\prime} ( t ) Y ( t ) - g ( t )
                  \frac{\D d}{\D dt} Y ( t ) }{\D g^{\,\prime} ( t ) } \right ) $ \\

for $t \in ( 0, \infty )$. Indeed, in view of (4.43), $g^{\,\prime} ( t )$ may be singular at
$t = 0$. \\
Finally, in view of (4.46), the solution $Y$ of the ODE in (4.49) does satisfy the initial
condition in (4.40).

\hfill $\Box$ \\

It is easy to see, Rosinger [1, pp. 199,200], that a lot more general examples than above in
(4.41) - (4.49), (4.40) can be constructed for {\it arbitrary} $f \in {\cal C}^\infty (
\mathbb{R}^l, \mathbb{R}^l ) \setminus {\cal D}iff^\infty ( \mathbb{R}^l, \mathbb{R}^l )$. \\

The conclusion is that, as indicated in (4.40), a class of {\it implicit non-autonomous} ODEs
which are {\it singular} at $t = 0$, is needed in order to be able, through the standard
construction in section 5, to obtain genuine Lie semigroups, or GLS-s. \\

{\bf Nonremovable Singularities} \\

As seen in section 5, the property of the actions $H$, see (4.30), and in particular (4.29),
(4.45) \\

(4.50) $~~~ H \in {\cal C}^0 ( [ 0, \infty ) \times M, M ) \setminus
                                          {\cal C}^1 ( [ 0, \infty ) \times M, M ) $ \\

is needed in order to obtain genuine Lie semigroups, or GLS-s. This property, as seen above,
corresponds to a singularity at $t = 0$ of the non-autonomous ODEs satisfied by such actions
$H$. Let us further note that such singularities of the actions $H$ must, therefore, be {\it
nonremovable}, and the example below shows that this is possible to attain. \\

Returning to the action $H$ in (4.29) with $M = \mathbb{R}$, we can associate with it the
action $K \in {\cal C}^\infty ( \mathbb{R} \times M, M )$, given by \\

(4.51) $~~~ K ( s, y ) ~=~ y + s y^2,~~~ s \in \mathbb{R},~~ y \in M $ \\

and then we have \\

(4.52) $~~~ H ( s^2, ~.~ ) ~=~ K ( | s |, ~.~ ),~~~ H ( t, ~.~ ) ~=~ K ( \sqrt t, ~.~ ) $ \\

for $s \in \mathbb{R},~ t \in [ 0, \infty )$. \\

However, even if $K$ has no singularity for $( s, y ) \in \mathbb{R} \times M$, the fact
remains that in the sense of $* * * )$ in (4.39), $H$ is singular at $t = 0$, as it fails to
be ${\cal C}^1$-smooth in a neighbourhood of $t = 0 $. \\ \\

{\bf 5. Standard Reduction to Autonomous ODEs and Evolution} \\

Here we present the needed details on the standard way mentioned briefly in Remark 4.1., and
according to which non-autonomous ODEs can be reduced to autonomous ones, with the all
important corresponding reduction of non-autonomous evolutions to an autonomous ones. This
standard reduction is fundamental in obtaining the general result in ( dim + 1 ) in section 4
above, which is a stepping stone towards the ultimate result in ( GLS ), a result still
open. \\

{\bf Autonomous Evolution.} Let be given the {\it autonomous} nonlinear system of ODEs \\

(5.1) $~~~ \begin{array}{l}
                \frac{\D d}{\D d t} Y ( t ) ~=~ F ( Y ( t ) ),~~~ t \in \mathbb{R} \\ \\
                Y ( t _ 0 ) ~=~ y_0
           \end{array} $ \\

with $F \in {\cal C}^\infty ( \mathbb{R}^l, \mathbb{R}^l ),~ t_0 \in \mathbb{R},~ y_0 \in
\mathbb{R}^l$, where the sought after solution is $Y \in {\cal C}^\infty ( \mathbb{R},
\mathbb{R}^l )$. For convenience, we assume that the unique solution $Y$ exists globally on
$\mathbb{R}$ for every initial condition $y_0 \in \mathbb{R}^l$. Thus we can associate with
(5.1) the {\it evolution operator} \\

(5.2) $~~~ \begin{array}{l}
               \mathbb{R} \ni t ~~\longmapsto~~~ E ( t ) : \mathbb{R}^l
                                         ~~\longrightarrow~~ \mathbb{R}^l \\ \\
               E ( t - t_0 ) ( y ) ~=~ Y ( t ),~~~
                                 t_0, t \in \mathbb{R},~~ y \in \mathbb{R}^l
            \end{array} $ \\

which, as mentioned, defines the {\it one dimensional Lie group action} on $\mathbb{R}^l$ \\

(5.3) $~~~ \mathbb{R} \times \mathbb{R}^l \ni ( t, y )
                   ~~\longmapsto~~ E ( t ) ( y ) \in \mathbb{R}^l $ \\

Thus the evolution operator $E$ has the {\it group properties} \\

(5.4) $~~~ \begin{array}{l}
                     E(0) = id_{\mathbb{R}^l} \\
                     \\
                     E(t + s) = E(t) E(s),~~~
                       t, s \in \mathbb{R}
                 \end{array} $ \\

in other words, we have the group homomorphism \\

(5.5) $~~~ \mathbb{R} \ni t ~\longmapsto~ E(t) \in
                                         {\cal D}iff^\infty (\mathbb{R}^l) $ \\

So much for a recapitulation of needed well known properties of autonomous nonlinear systems
of ODEs. \\

{\bf Non-autonomous Evolution.} And now, let us consider the {\it non-autonomous} nonlinear
system of ODEs \\

(5.6) $~~~ \begin{array}{l}
                    \frac{\D d}{\D d t} Y(t) =
                       F(t,Y(t)),~~ t \in \mathbb{R} \\
                    \\
                    Y(t_0) = y_0 \in \mathbb{R}^l
                 \end{array} $ \\

where $F \in {\cal C}^\infty(\mathbb{R}^{l+1},\mathbb{R}^l),~ t_0 \in \mathbb{R},~ y_0 \in
\mathbb{R}^l$, and for convenience, the unique solution $Y \in {\cal C}^\infty(\mathbb{R},
\mathbb{R}^l)$ is supposed to exist for all $t \in \mathbb{R}$. This time, the associated {\it
non-autonomous evolution operator} $E \in {\cal C}^\infty( \mathbb{R}^2 \times \mathbb{R}^l,
\mathbb{R}^l)$ has the property \\

(5.7) $~~~ \begin{array}{l}
                     \mathbb{R} \times \mathbb{R} \ni (t_0,t)
                        ~\longmapsto~ E(t_0,t) : \mathbb{R}^l
                                             \longrightarrow \mathbb{R}^l   \\
                     \\
                     E(t_0,t)(y_0) = Y(t),~~~
                        t_0, t \in \mathbb{R},~ y_0 \in \mathbb{R}^l
                 \end{array} $ \\

However, it is important to note that if we solve only for a given fixed $t_0 \in \mathbb{R}$
the non-autonomous ODE system in (5.6) and with all the initial conditions $y_0 \in
\mathbb{R}^l$, then (5.7) will in general {\it not} give the full information on $E$, and even
less will give it in an explicit manner. This is unlike the case with autonomous ODEs, see
(5.2). \\

Nevertheless, as seen in (5.22), the full information on the non-autonomous evolution operator
$E$ can be recovered from the knowledge of (5.7) even for one single $t_0 \in \mathbb{R}$, and
moreover, it can be recovered without having to solve the ODE system (5.6) for every other
$t_0 \in \mathbb{R}$, but only by solving some additional algebraic equations. Indeed, this
becomes possible, once we transform the non-autonomous ODE system in (5.6) into the autonomous
ODE system (5.10), and then we use the corresponding autonomous evolution operator $E_A$. \\

In the case of (5.6), (5.7), the relations (5.4) take the following more general, {\it
non-autonomous form of group property}, namely \\

(5.8) $~~~ \begin{array}{l}
                    E(t,t) = id_{\mathbb{R}^l},~~ t \in \mathbb{R} \\ \\
                    E(s,r) E(t,s) = E(t,r),~
                    t, s, r \in \mathbb{R}
                 \end{array} $ \\ \\

{\bf Reduction to Autonomous Evolution.} Now we recall the standard way the non-autonomous
system (5.6) can be {\it reduced} to an autonomous one, such as for instance in (5.1). \\

Namely, let us augment the function $Y \in {\cal C}^\infty(\mathbb{R},\mathbb{R}^l)$ in (5.6)
to the function $Y_A \in {\cal C}^\infty(\mathbb{R},\mathbb{R}^{l+1})$, given by \\

(5.9) $~~~ Y_A(t) = (t,Y(t)),~~~ t \in \mathbb{R} $ \\

then clearly, the $l$-dimensional non-autonomous ODE system (5.6) is {\it equivalent} with the
$l+1$-dimensional autonomous ODE system \\

(5.10) $~~~ \begin{array}{l}
                   \frac{\D d}{\D d t} Y_A(t) =
                       (1, F(Y_A(t))),~~~ t \in \mathbb{R} \\
                   \\
                   Y_A(t_0) = (t_0,y_0) \in \mathbb{R}^{l+1}
                \end{array} $ \\

Two facts should be noted here. \\

First, the {\it non-autonomous infinitesimal generator} $F$ in (5.6), becomes associated by
(5.10) with the {\it autonomous infinitesimal generator} $F_A \in {\cal C}^\infty(
\mathbb{R}^{l+1},\mathbb{R}^{l+1})$ defined by \\

(5.11) $~~~ F_A(y_A) = (1, F(y_A)),~~~ y_A \in \mathbb{R}^{l+1} $ \\

Second, in the initial condition $(t_0,y_0)$ in (5.10), the first coordinate $t_0 \in
\mathbb{R}$ must be the {\it same} with that in the left hand term. Clearly, this is a direct
consequence of the definition of $Y_A$ in (5.9). \\

Let now $E_A \in {\cal C}^\infty(\mathbb{R} \times \mathbb{R}^{l+1},\mathbb{R}^{l+1})$
be the autonomous evolution operator associated with (5.10), then \\

(5.12) $~~~ E_A(t - t_0)(t_0,y_0) =
               Y_A(t),~~~ t_0, t \in \mathbb{R},~ y_0 \in \mathbb{R}^l $ \\

and it is easy to see that it can be decomposed as follows \\

(5.13) $~~~ E_A(s)(t,y) = (E_{A 1}(s)(t,y),
                    E_{A 2}(s)(t,y)),~~ t, s \in \mathbb{R},~ y \in \mathbb{R}^l $ \\

with the ${\cal C}^\infty$-smooth functions \\

(5.14) $~~~ \begin{array}{l}
                      \mathbb{R} \ni s ~\longmapsto~ E_{A 1}(s) :
                                  \mathbb{R}^{l+1} \longrightarrow \mathbb{R} \\
                      \\
                      \mathbb{R} \ni s ~\longmapsto~ E_{A 2}(s) :
                                  \mathbb{R}^{l+1} \longrightarrow \mathbb{R}^l
                  \end{array} $ \\

In view of (5.7), (5.9), (5.12), (5.13), it follows that, for $t_0, t \in \mathbb{R},~ y_0 \in
\mathbb{R}^l$, we have \\

(5.15) $~~~ E(t_0,t)(y_0) = Y(t) = E_{A 2}(t - t_0)(t_0,y_0) $ \\

Also, $E_A$ satisfies the group properties corresponding to the autonomous case (5.4),
namely \\

(5.16) $~~~ \begin{array}{l}
                     E_A(0) = id_{\mathbb{R}^{l+1}} \\
                     \\
                     E_A(t + s) = E_A(s) E_A(t),~~~ t, s \in \mathbb{R}
            \end{array} $ \\

which result in the group homomorphism \\

(5.17) $~~~ \mathbb{R} \ni t ~\longmapsto~ E_A(t) \in
                                       {\cal D}iff^\infty(\mathbb{R}^{l+1}) $ \\

Therefore, for $t, s, r \in \mathbb{R}$ and $y \in \mathbb{R}^l$, we obtain the relations \\

(5.18) $~~~ \begin{array}{l}
                      E_{A 1}(r + s)(t,y) =
                         E_{A 1}(r)(E_{A 1}(s)(t,y),
                                           E_{A 2}(s)(t,y)) \\
                      \\
                      E_{A 2}(r + s)(t,y) =
                         E_{A 2}(r)(E_{A 1}(s)(t,y),
                                            E_{A 2}(s)(t,y))
            \end{array} $ \\

And in view of (5.7) - (5.12), we obtain for $t, s \in \mathbb{R},~ y \in \mathbb{R}^l$ \\

(5.19) $~~ \begin{array}{l}
                     E_{A 1}(s)(t,y) = t + s \\
                     \\
                     E_{A 2}(s)(t,y) =
                               E(t,t + s)(y)
           \end{array} $ \\

In particular, we can check for $E$ the non-autonomous form of the semigroup property (5.8).
Namely, for $t, s, r \in \mathbb{R},~ y \in \mathbb{R}^l$, we have \\

(5.20) $~~~ \begin{array}{l}
                   E(t,t + s + r)(y) = \\
                   ~~~~~~~~~~~~~~~~~~~= E_{A 2}(s + r)(t,y) =
                   \\
                   ~~~~~~~~~~~~~~~~~~~= E_{A 2}(r)(t + s,
                                 E_{A 2}(s)(t,y)) = \\
                   ~~~~~~~~~~~~~~~~~~~= E_{A 2}(r)(t + s,
                                 E(t, t + s)(y)) =  \\
                   = E(t + s,t + s + r) E(t, t + s)(y)
           \end{array} $ \\

{\bf Recovering the Non-autonomous Evolution.} Let us show now how we can recover the full
non-autonomous evolution operator $E$ in (5.8), from anyone of its particular cases in (5.7),
which corresponds merely to a certain $t_0 \in \mathbb{R}$ fixed. For that purpose, we shall
use the autonomous extension $E_A$ of $E$. \\

Indeed, (5.15) gives \\

(5.21) $~~~ E(t,s)(y) = E_{A 2}(s,t)(t,y),
                          ~~~ t,~ s \in \mathbb{R},~ y \in \mathbb{R}^l $ \\

We fix now $t_0 \in \mathbb{R}$. Given any $t \in \mathbb{R}$, if we can compute the
mapping \\

$~~~~~~ \mathbb{R}^l \ni y \longmapsto y_* \in \mathbb{R}^l $ \\

where $y_*$ is the solution of the algebraic equation \\

$~~~~~~ E(t_0,t)(y_*) = y $ \\

then it is clear that we thus obtain $(t,y) = E_A(t - t_0)(t_0,y_0)$, hence \\

$~~~~~~ E(t,s)(y) = E_{A 2}(s - t) E_A(t - t_0)(t_0,y_*) = E(t_0,s)(y_*) $ \\

Therefore \\

(5.22) $~~~ \begin{array}{l}
                       E(t,s)(y) = E(t_0,s)(y_*),~~~
                         t, s \in \mathbb{R},~ y \in \mathbb{R}^l    \\
                       \\
                       \mbox{where}~~ y_* ~~\mbox{is a solution of}~~ E(t_0,t)(y_*) = y
            \end{array} $ \\ \\

We shall also assume that \\

(5.23) \quad $ y_* ~~\mbox{depends}~~ {\cal C}^\infty-\mbox{smoothly on}~~ t,~ y $ \\

Let us illustrate the above with a simple example. We consider the non-autonomous ODE \\

$~~~~~~ \begin{array}{l}
         \frac{\D d}{\D d t} Y(t) = 2 t,~~ t \in \mathbb{R} \\
         \\
         Y(t_0) = y_0 \in \mathbb{R}
        \end{array} $ \\

for which obviously \\

$~~~~~~ \begin{array}{l}
                Y(t) = E(t_0,t)(y_0) = E_{A 2}(t - t_0)(t_0,y_0) = t^2 - t_0^2 + y_0  \\ \\
                Y_A(t) = (t,E_{A 2}(t - t_0)(t_0, y_0)) \\ \\
                E_{A 1}(s)(t,y) = t + s \\ \\
                E_{A 2}(s)(t,y) = E(t,t + s)(y) = s^2 + 2  s t + y
        \end{array} $ \\

We now check the group properties (5.15). First, we note that the relation \\

$~~~~~~ E_A(0) = id_{\mathbb{R}^{l+1}} $ \\

results immediately from the expressions of $E_{A 1}$ and $E_{A 2}$ above. Then \\

$~~~~~~ E_A(r) E_A(s) = E_A(s + r),~~~ s, r \in \mathbb{R} $ \\

holds since the previous relations give $ E_A(s)(t,y) = (t + s,s^2 + 2 s t + y)$, hence
$E_A(r) E_A(s)(t,y) = E_A(r)(t + s,s^2 + 2 s t + y) = (t + s + r,r^2 + 2 r (t + s) + s^2 +
2 s t + y)$, while $E_A(s + r)(t,y) = (t + s + r, (s + r)^2 + 2(s +r) t + y)$. \\

In other words, the non-autonomous version (5.8) of the group property works as follows, for
$t, s, r \geq 0,~ y \in M$, we have \\

$~~~~~~ E(t,s)(y) = s^2 - t^2 + y $ \\

hence clearly $E(t,t)(y) = y$, while \\

$~~~~~~ E(s,r) E(t,s)(y) = r^2 - s^2 + (s^2 - t^2 + y) = r^2 - t^2 + y = E(t,r)(y) $ \\

The autonomous version, which corresponds to the non-autonomous one, according to (5.9) -
(5.22), will have the group property (5.16) working as follows. For $t, s, r \geq 0,~ y \in
M$, we have \\

$~~~~~~ E_A(s)(t,y) = (t + s,s^2 + 2 s t + y) $ \\

thus obviously $E_A(0)(t,y) = (t,y)$, while \\

$~~~~~~ \begin{array}{l}
            E_A(r) E_A(s)(t,y) = (t + s + r, r^2 + 2 r (t + s) + s^2 + 2 s t + y) = \\ \\
            ~~~~~~~~~~~~ = (t + s + r,(s + r)^2 + 2(s + r) t + y) = E_A(s + r)(t,y)
        \end{array} $ \\

Finally, in the case of this example, the relation (5.22) works as follows. Given a fixed
$t_0 \in \mathbb{R}$, the algebraic equation \\

$~~~~~~ t^2 - t_0^2 + y_* = y $ \\

can obviously be solved in $y_*$ for every $t, y \in \mathbb{R}$, and it gives $y_* = y - t^2
+ t_0^2$, which also satisfies (5.23), therefore \\

$~~~~~~ E(t,s)(y) = E(t_0,s)(y_*) = s^2 - t_0^2 + y_* = s^2 - t_0^2 + y $ \\ \\

{\bf 6. Examples of Genuine Lie Semigroups} \\

We return to the examples in section 4, and show the way in which the respective actions can
be associated with {\it genuine Lie semigroup} actions, by applying to them the method of
reduction in section 5. \\

First we note that the results in section 5, where the ODEs are defined on the whole of
$\mathbb{R}$, have to be adapted since, in general, and as seen with the examples in section 4,
one no longer has $t \in \mathbb{R}$, but only $t \in (0,\infty)$, or at most $t \in [0,
\infty)$, see (4.29), (4.43) - (4.46), or Rosinger [1, p. 199, (13.2.36) - (13.2.38)]. As
noted in section 4, this restriction of the domain of $t$ is due to the {\it nonremovable
singularities} of the respective ODEs, and more specifically, of their solutions of
interest. \\

Concerning the autonomous case (5.1) - (5.5), such a restriction of the domain of $t$ means
that the respective evolution operator $E$ in (5.2) will only be defined for $t \in [0,
\infty)$, and instead of the group property (5.4), will only have the {\it semigroup}
property \\

(6.1) $~~~ E(t + s) = E(s) E(t),~~~ t,~ s \in [0,\infty) $ \\

In particular, $E(t)$, with $t \in (0,\infty)$, may fail to be invertible, since $E(- t)$ need
not exist, and thus, we could not always obtain from (5.4) the relations $E(t) E(- t) =
E(- t) E(t) = id_{\mathbb{R}^l}$. \\

In the non-autonomous case (5.6) - (5.8), the evolution operator $E$ in (5.7) will only be
defined for $t_0,~ t \in [0,\infty)$, and will have the following {\it semigroup} version of
property (5.8) \\

(6.2) $~~~ E(r,s) E(t,s) = E(t,r),~~~ t,~ s,~ r \in [0,\infty) $ \\

This however allows for the existence of its {\it inverses}, since we have for any $t,~ s \in
[0,\infty)$ \\

(6.3) $~~~~ \begin{array}{l}
                   E(t,s) E(s,t) = E(s,s) = id_{\mathbb{R}^l} \\ \\
                   E(s,t) E(t,s) = E(t,t) = id_{\mathbb{R}^l}
            \end{array} $ \\

Let us now transform, more precisely reduce, the thus restricted non-autonomous version of
(5.6) - (5.8) into the corresponding restricted autonomous version of (5.9) - (5.22). Then
clearly, for $t,~ s \in [0,\infty),~ y \in \mathbb{R}^l$, we obtain \\

(6.4) $~~~ E_A(s)(t,y) = (t + s, E(t, t + s)(y)) $ \\

or equivalently \\

(6.5) $~~~ E_A(s - t)(t,y) = (s, E(t, s)(y)),~~~ 0 \leq t \leq s < \infty $ \\

Therefore, in the case of such singular non-autonomous ODEs and of their corresponding
evolution operators $E$, after the transformation into the autonomous case, except for the
trivial situation of $s = t$, we need no longer be able to benefit from the existence of the
inverses in (6.3), when we deal with the associated autonomous evolution operator $E_A$. \\

This is precisely at the basis of our construction of {\it genuine Lie semigroup} actions. \\

And in particular, this is how we shall associate such genuine Lie semigroup actions with the
examples in section 4. \\

Before going further, let us note that in the case of both non-autonomous examples in section
4, we do have the corresponding versions of (6.2), namely \\

(6.6) $~~~ E(s,r) E(t,s) = E(t,r),~~~ t,~ s,~ r \in [0,\infty) $ \\

Indeed, for the example in (4.29) - (4.36), this follows from the fact that $H$ in (4.29) is a
solution on $(0,\infty)$ of (4.31) - (4.36), while in addition, see (4.37), $H$ is such that
$H(0,y) = y$, for all $y \in M$, and also $H \in {\cal C}^0([0,\infty) \times M,M)$.

A similar argument will apply to the more general examples in (4.43) - (4.46), or Rosinger [1,
p. 199, (13.2.36) - (13.2.38)]. \\

In view of (6.6), we obtain for the example in (4.29) - (4.36), the relation \\

(6.7) $~~~ E(0,t)(y) = y + \sqrt{t} y^2,~~~ t \in [0,\infty),~ y \in M $ \\

while for the example in (4.43) - (4.49), we shall have \\

(6.8) $~~~ E(0,t)(y) = ( 1 - g(t) ) y + g(t) f(y),~~~ t \in [0,\infty),~
                                                                  y \in \mathbb{R}^l $ \\

These two relations will help us in fully computing the respective non-autonomous evolution
operators $E$ for the mentioned examples. \\

Indeed, let us determine $E(t,s)(y)$ for the first example, and do so for all $t,~ s \in (0,
\infty),~ y \in M = \mathbb{R}$. From (6.6) we have \\

$~~~~~~ E(t,s)(y) = E(0,s) E(t,0)(y) = E(0,s) ( E(0,t) )^{-1} (y) $ \\

Hence, proceeding for $t_0 = 0$ as in (5.22), let us assume that \\

$~~~~~~ ( E(0,t) )^{-1} (y) = y_* \in M $ \\

then clearly \\

$~~~~~~ y = E(0,t)(y_*) = y_* + \sqrt{t} y_*^2 $ \\

and we note the important consequence that \\

$~~~~~~ \lim_{~t \to 0}~y_* = y \in M $ \\

provided that $y_*$ is bounded. \\

Now computing $y_*$ from the above quadratic equation, we obtain \\

$~~~~~~ y_* = \frac{\D -1 \pm \sqrt{1+4 \sqrt{t} y}}{\D 2 \sqrt{t}}\,,~~~
                       t \in (0,\infty),~ y \in M,~ 1+4 \sqrt{t} y \geq 0 $ \\

and in order to secure the above limit, it follows that we must choose \\

$~~~~~~ \begin{array}{l}
               y_* = \frac{\D -1 + \sqrt{1+4 \sqrt{t} y}}{\D 2 \sqrt{t}} =
                  \frac{\D 2 y}{\D 1+ \sqrt{1+4 \sqrt{t} y}}\,, \\ \\
    ~~~~~~~~~~~~~~~~~~~~~~~~~~~~~~~~~~ t \in (0,\infty),~ y \in M,~ 1+4 \sqrt{t} y \geq 0
         \end{array} $ \\

In this way (5.23) is satisfied, while

\bigskip
(6.9) \quad $ \begin{array}{l}
                  E(t,s)(y) = E(0,s)(y_*) = \\
                  \\
                  ~~~~~~~~~~
                  = \frac{\D 2 y}{\D 1+ \sqrt{1+4 \sqrt{t} y}} +
                  \sqrt{s}~ \frac{\D 4 y^2}{\D
                  \left ( 1+ \sqrt{1+4 \sqrt{t} y} \right ) ^2} \\
                  \\
                  ~~~~~~~~~~~~~~~~~~~~~~~~~~~~~~~~~~~~~~
                         t,~ s \in (0,\infty),~ y \in M,~
                                          1+4 \sqrt{t} y \geq 0
                  \end{array} $ \\

and we have obtained the {\it full expression of the non-autonomous evolution operator} $E$
for the first example (4.29) - (4.36) in section 4. \\

The {\it genuine Lie semigroup} of actions generated by this example in (4.29) - (4.36) will
now be given by the {\it autonomous evolution operator} $E_A$, which corresponds to $E$ above,
according to (5.9), namely, see also (6.4) \\

(6.10) $~~~ \begin{array}{l}
                     E_A(s)(t,y) = (t + s, E(t, t + s)(y) ) \\ \\
                     ~~~~~~~~~~~~~~~~~~~~~~~~~~~~ t,~ s \in [0,\infty),~ y \in M,~ 1+4
                                               \sqrt{t} y \geq 0
            \end{array} $ \\

The genuine Lie semigroup actions generated by the example in (4.43) - (4.49) can be obtained
in a similar manner, provided that one can solve in $y_0 \in \mathbb{R}^l$ the corresponding
algebraic equations \\

$~~~~~~ y = ( 1 - g(t) ) y_0 + g(t) f(y_0) $ \\

for given $(t,y) \in (0,\infty) \times \mathbb{R}^l$. \\

{\bf Milder Singularities.} We show with an example that the condition \\

( SING ) $~~~ H ( t, ~.~ ) \notin {\cal D}iff^\infty ( M ),~~~ t \in ( 0, \infty ) $ \\

is {\it not} necessary, in order to obtain genuine Lie group actions by using the above method.
Indeed, let $M = \mathbb{R}$ and $f \in {\cal C}^\infty(M,M)$ be given by \\

(6.11) $~~~ f(y) = \frac{\D 1}{\D y^2 + 1},~~~ y \in M $ \\

We define $H \in {\cal C}^0([0,\infty) \times M,M) \cap {\cal C}^\infty ((0,\infty) \times M,
M )$ such that \\

(6.12) $~~~ H(t,y) = ( 1 - \sqrt{t} ) y + \sqrt{t} f(y) =
                  ( 1 - \sqrt{t} ) y + \frac{\D \sqrt{t}}{\D y^2 + 1},~~~ \\ \\
 ~~~~~~~~~~~~~~~~~~~~~~~~~~~~~~~~~~~~~~~~~~~~~~~~~~~~~~~~~~~~~t \in [0,\infty), ~ y \in M $ \\

Then it follows that \\

(6.13) $~~~ \begin{array}{l}
                     H(1,~.~) = f \notin {\cal D}iff^\infty(M) \\ \\
                     H(t,~.~) \in {\cal D}iff^\infty(M)
                     ~~~\mbox{for}~~ t \in [0,4/9) \cup (4,\infty)
             \end{array} $ \\

However, it is obvious that the above procedure applied to the examples in section 4, is
equally applicable to (6.11) - (6.13), and again, it will lead to genuine Lie semigroup
actions. \\ \\

{\bf 7. Singularity, Continuity, Smoothness and Domains \\
\hspace*{0.4cm} of Action} \\

In view of sections 4 - 6, one possible method to obtain genuine Lie semigroups is that given
by the evolution operators $E_A$ of the autonomous singular ODEs, which are associated in the
above standard manner with the non-autonomous singular ODEs in (4.37) - (4.40). \\
In other words, such genuine Lie semigroup actions on suitable subsets $\widetilde{M}$ in
Euclidean spaces are given by mappings \\

(7.1) $~~~ E_A \in ({\cal C}^0([0,\infty) \times \widetilde{M},\widetilde{M}) \cap
                          {\cal C}^\infty((0,\infty) \times \widetilde{M},\widetilde{M}))
                          \setminus  \\
~~~~~~~~~~~~~~~~~~~~~~~~~~~~~~~~~ {\cal C}^1([0,\infty) \times \widetilde{M}, \widetilde{M}) $ \\

With respect to the domains of action of such genuine Lie semigroups, as constructed in
sections 4 - 6, we have to note the following. We have started with certain open subsets $M$
in Euclidean spaces, and see for instance (4.29), with singular actions \\

(7.2) $~~~ H \in ({\cal C}^0([0,\infty) \times M,M) \cap {\cal C}^
                        \infty((0,\infty) \times M,M)) \setminus \\
~~~~~~~~~~~~~~~~~~~~~~~~~~~~~~~~~~ {\cal C}^1([0,\infty) \times M,M) $ \\

which clearly were not any kind of semigroup actions on the respective open subsets $M$. \\

Then, we found non-autonomous singular ODEs which were satisfied by these singular actions $H$,
and associated with them in the standard manner autonomous singular ODEs. \\

Finally, the evolution operators $E_A$ of these associated autonomous singular ODEs gave us
the genuine Lie semigroup actions. \\

However, such an $E_A$ is {\it no longer} acting on $M$, but on the set with one dimension
higher, namely \\

(7.3) $~~~ \widetilde{M} = [0,\infty) \times M $ \\

Thus we are led to the {\bf Open Problem} ( GLS ) formulated at the beginning of section
4. \\ \\

{\bf 8. Remark on Singularities} \\

Let us note that in order to obtain genuine Lie semigroups, we can use {\it milder} forms of
singularities than those in section 4. For instance, instead of $H$ given by (4.29), let us
consider it defined as follows \\

(8.1) \quad $ H(t,y) = y + t y^2,~~~ t \in \mathbb{R},~~~ y \in M = \mathbb{R} $ \\

Then (4.30) becomes replaced with \\

(8.2) \quad $ \begin{array}{l}
                   H(0,y) = y,~~~ y \in M \\
                   \\
                   H(t,~.~) \notin {\cal D}iff^\infty(M),~~~
                                t \in \mathbb{R} \setminus \{ 0 \} \\
                   \\
                   H \in {\cal C}^\infty(\mathbb{R} \times M,M)
                \end{array} $ \\

while the action in (4.29) limited to $t \in [ 0, \infty )$, extends now to the following one
defined for all $t \in \mathbb{R}$, namely \\

(8.3) \quad $ \mathbb{R} \times M \ni (t,y) \longmapsto H(t,y) \in M $ \\

However, this extended action still {\it cannot} be part of a group or local group action on
$M$, since in view of (8.2), $H(t,~.~)$, with $t \in \mathbb{R} \setminus \{ 0 \}$, is not a
${\cal C}^\infty$-smooth diffeomorphism of $M$. \\

Proceeding now in a manner similar with that in (4.31) - (4.36), it follows that $H$ in (8.1),
as a function of $t \in \mathbb{R}$, and for any given fixed $y \in M$, will satisfy the ODE \\

(8.4) \quad $ \begin{array}{l}
                    \frac{\D \partial}{\D \partial t} Y(t) =
                    \frac{\D 2 Y(t)^2}{\D 1 + 2 t Y(t)
                    + \sqrt{1 + 4 t Y(t)}}\,, \\
                    \\
                    ~~~~~~~~~~~~~~~~~~~~~~~~~~~~~~~~~~~~~~~~~~
                    t \in \mathbb{R},~ 1 + 4 t Y(t) \geq 0
               \end{array} $ \\

if \\

(8.5) \quad $ 1 + 2 t y \geq 0 $ \\

while it will satisfy the ODE \\

(8.6) \quad $ \begin{array}{l}
                    \frac{\D \partial}{\D \partial t} Y(t) =
                    \frac{\D 1 + 2 t Y(t) + \sqrt{1 + 4
                    t Y(t)}}{\D 2 t^2}\,, \\
                    \\
                    ~~~~~~~~~~~~~~~~~~~~~~~~~~~~~~~~~~~~~~~~~~
                    t \in \mathbb{R} \setminus \{ 0 \},~ 1 + 4 t
                    Y(t) \geq 0
                \end{array} $ \\

if \\

(8.7) \quad $ 1 + 2 t y \leq 0 $ \\

Furthermore, both these ODEs will for $H$ specified above be associated with the initial
condition \\

(8.8) \quad $ \lim_{~t \to 0}~ Y(t) = y \in M $ \\

We can note in the above example (8.1) - (8.8) that in the case of (8.5), the corresponding
ODE in (8.4) which is satisfied by $H$, seen as a function of $t \in \mathbb{R}$, and with
$y \in M$ fixed, is {\it not} singular at $t = 0$.

Nevertheless, the respective solution $H$, for $t \in \mathbb{R} \setminus \{ 0 \}$, still
{\it cannot} be part of a group or local group of transformations on $M$, in view of (
8.2). \\ \\

{\bf 9. Evolution PDEs and Genuine Lie Semigroups} \\

In this section, as a possible alternative to the method in section 4, we give an indication
about another way one may try to generate genuine Lie semigroup actions. This alternative way
is suggested by the well established literature on solving initial value problems for
evolution PDEs through the associated semigroups of operators acting on the respective initial
values. The early basic result in the case of linear evolution PDEs is the celebrated
Hille-Yoshida theorem, which was followed by a rather large body of more recent results,
including nonlinear developments, see for instance Pazy and the references cited there. \\

The important point to note here is that, in general, the solutions of the initial value
problems for evolution PDEs will be given by {\it semigroups}, rather than groups, of such
operators. A well known class of evolution PDEs for which, typically, one can only obtain such
semigroups of operators is that of parabolic equations. \\

The idea suggested in this section is to try to use such semigroups acting on initial values,
in order to generate the genuine Lie semigroups which are the object of our interest in this
paper. \\

For simplicity, we shall again consider the case when $M = \mathbb{R}^l$ is the open set on
which we want to define a genuine Lie semigroup action. \\

Let us therefore take any evolution PDE of the form \\

(9.1) \quad $ D_t U(t,x) = T(x,D) U(t,x),~~~ t \in [0,\infty),~ x \in M $ \\

where $U \in {\cal C}^\infty([0,\infty) \times M, \mathbb{R})$ is the unknown function, while
$T(x,D)$ is a partial differential operator in $x$ alone. \\

Associated with (9.1), we consider the initial value problem \\

(9.2) \quad $ U(0,x) = f(x),~~~ x \in M $ \\

where $f$ belongs to a suitable class of functions, namely \\

(9.3) \quad $ f \in {\cal F}(M) \subseteq {\cal C}^\infty(M,\mathbb{R}) $ \\

We shall assume that there exists a semigroup of operators acting on the typically {\it
infinite} dimensional vector space ${\cal F}(M)$, namely \\

(9.4) \quad $ [0,\infty) \ni t \longmapsto E(t) : {\cal F}(M) \rightarrow {\cal F}(M) $ \\

such that, given any $f \in {\cal F}(M)$, if we define \\

(9.5) \quad $ U(t,x) = (E(t) f)(x),~~~ t \in [0,\infty),~ x \in M $ \\

then $U$ is a solution of (9.1), (9.2). \\

So far, we have been moving within the well established framework of the mentioned literature,
provided that we work with suitable evolution equations (9.1) and initial values (9.2). \\

Now the idea which is the subject of this section is to make the following further assumption.
Namely, suppose that there exists a family of functions \\

(9.6) \quad $ V_{a,b} \in {\cal F}(M),~~~ a \in M,~ b \in B $ \\

where $B \subseteq \mathbb{R}^k$ is a suitable open subset, such that \\

(9.7) \quad $ E(t) V_{a,b} = V_{\alpha(t,a,b),\beta(t,a,b)},~~~
                                    t \in [0,\infty),~ a \in M,~ b \in B $ \\

for appropriate functions \\

(9.8) \quad $ \alpha \in {\cal C}^\infty([0,\infty) \times M \times B,
                      M),~ \beta \in {\cal C}^\infty([0,\infty) \times M \times B,B) $ \\

In this case, the assumed semigroup property \\

(9.9) \quad $ E(s) E(t) = E(t + s),~~~ t, s \in [0,\infty) $ \\

together with (9.6), (9.7) will result in \\

(9.10) \quad $ \begin{array}{l}
                   \alpha(t+s,a,b) = \alpha(s,\alpha(t,a,b),\beta(t,a,b)) \\
                   \\
                   \beta(t+s,a,b) = \beta(s, \alpha(t,a,b),\beta(t,a,b))
                \end{array} $ \\

for $t,s \in [0,\infty),~ a \in M,~ b \in B$. \\

It follows that if, for instance \\

(9.11) \quad $ \beta(t,a,b) = b,~~~ t \in [0,\infty),~ a \in M,~ b \in B $ \\

then (9.10) yields \\

(9.12) \quad $ \alpha(t+s,a,b) = \alpha(s,\alpha(t,a,b),b),~~ t,s \in [0,\infty),~ a \in M,~ b \in B $ \\

hence for any fixed $b \in B$, the mapping \\

(9.13) \quad $ [0,\infty) \times M \ni (t,a) \longmapsto \alpha(t,a,b) \in M $ \\

is a semigroup action on $M$, and it is a priori {\it not} impossible that it may indeed be a
{\it genuine Lie semigroup} action. \\

It is easy to see that the assumptions (9.6) - (9.8) can be satisfied by a large variety of
{\it soliton} solutions, for instance. Indeed, for simplicity, let us consider the one
dimensional case, when $M = \mathbb{R}$, and assume that (9.1) has a soliton solution \\

(9.14) \quad $ U(t,x) = W(x - c t),~~~ t \in [0,\infty),~ x \in M $ \\

for $c \in \mathbb{R}$. In this case (9.5) becomes \\

(9.15) \quad $ (E(t) W)(x) = W(x - c t),~~~ t \in [0,\infty),~ x \in M $ \\

For further more detailed illustration, let us consider the well known Burgers equation \\

(9.16) \quad $ U_t(t,x) + U(t,x) U_x(t,x) = \mu U_{xx}(t,x),~~~ t,x \in \mathbb{R} $ \\

where $\mu \in (0,\infty)$. This equation has a soliton solution given by \\

(9.17) \quad $ U(t,x) = c - \sqrt{c^2 + d}~ \tanh \left ( \frac{\D \sqrt{
                              c^2 + d}}{\D 2 \mu}~ (x - x_0 - c t) \right ), \\ \\
~~~~~~~~~~~~~~~~~~~~~~~~~~~~~~~~~~~~~~~~~~~~~~~~~~~~~~~~~~~~~~~~~~~~~~ t,x \in \mathbb{R} $ \\

where $x_0, c, d \in \mathbb{R}$ are arbitrary fixed, and $c^2 + d > 0$. For this soliton it
is clear that (9.6) - (9.8), (9.11) - (9.13) are satisfied, if we take \\

(9.18) \quad $ a = x_0 \in M = \mathbb{R},~~ b = (c,d) \in B \subset \mathbb{R}^2 $ \\

where $B = \{ (c,d)\in \mathbb{R}^2 ~|~ c^2 + d > 0 \}$, in which case \\

(9.19) \quad $ \begin{array}{l}
                     \alpha(t,a,b) = x_0 + c t \\
                     \\
                     \beta(t,a,b) = b
                  \end{array} $ \\

for $t \in [0,\infty),~ a \in M,~ b \in B$. \\ \\

{\bf 10. Other Instances of Semigroups of Actions} \\

In a private correspondence, P J Olver mentioned further instances in which semigroups of
actions appear in a natural way. Not all of them, however, need be genuine Lie semigroup
actions. \\

A first such example happens in the framework of (1.1), when a certain subset $S \subset M$ is
given, and we are interested in the set of actions which invariate it, namely \\

$~~~~~~ G_S = \{~ g \in G ~|~ g S \subseteq S ~\} $ \\

For instance, let \\

$~~~~~~ S = (-1,1) \times \mathbb{R} \subset \mathbb{R}^2 = M $ \\

and $G = (~(0,\infty),~.~)$ be the usual Lie group, which is supposed to act on $M$ according
to \\

$~~~~~~ G \times M \ni (g,(x,y)) \longmapsto (gx,y) \in M $ \\

Then clearly \\

$~~~~~~ G_S = (0,1] $ \\

which is a semigroup action. However, $G_S$ is not a genuine Lie semigroup action, since it is
a subsemigroup of the Lie group action $G$. \\

A second example is given by ODEs with inequality constraints. Let us, for instance, consider
the differential equation \\

$~~~~~~ \frac{\D d}{\D dx} U(x) = 0,~~~ x \in \mathbb{R} $ \\

with the inequality constraint \\

$~~~~~~ U(x) > 0,~~~ x \in \mathbb{R} $ \\

Then the Lie group actions \\

$~~~~~~ (x,u) \longmapsto (x,u+c) $ \\

are symmetries, only if $c \geq 0$. Needless to say, in view of a large class of applications,
such as control theory or differential games, for instance, where ODEs with inequality
constraints play a crucial role, the study of semigroups of symmetries of such equations can
present a special interest. \\

Another class of examples, this time related to PDEs, is given by generalized symmetries, see
Olver [1, chap. 5]. Indeed, the evolution PDEs governing the flow of a generalized symmetry
often only define a semigroup. A good example is the symmetry \\

$~~~~~ V = U_{xx} \partial _U $ \\

which is a symmetry of any linear constant coefficient PDE. Its flow is \\

$~~~~~~ U_t = U_{xx} $ \\ \\

{\bf Appendix} \\

In Rosinger [1] the {\it global} approach to arbitrary Lie group actions on smooth functions
was introduced and developed. And it was shown that for such a purpose, the use of a {\it
parametric} representation of functions upon which the Lie groups are supposed to act is
particularly appropriate. Here we present the essentials in this regard, as needed in this
paper. \\

Let us consider linear or nonlinear PDEs of the general form \\

(A.1) $~~~ T ( x, D )~ U ( x ) ~=~ 0,~~~ x \in \Omega $ \\

where $\Omega$ is an open subset in $\mathbb{R}^n$. \\

Lie group theory deals, among others, with those symmetries of solutions $U : \Omega
\longrightarrow \mathbb{R}$ of any given PDE in (A.1) which lead to other solutions of the
same equation. For that purpose, one takes $M = \Omega \times \mathbb{R}$ and finds the
corresponding Lie groups $G$ and their actions on $M$, namely \\

(A.2) $~~~ G \times M \ni ( g, ( x, u ) ) ~\longmapsto~ g ( x, u ) =
                                    ( g_1 ( x, u ), g_2 ( x, u ) ) \in M $ \\

where \\

(A.3) $~~~ \begin{array}{l}
                G \times M \ni ( g, ( x, u ) ) ~\longmapsto~ g_1 ( x, u ) \in \Omega \\
                G \times M \ni ( g, ( x, u ) ) ~\longmapsto~ g_2 ( x, u ) ) \in \mathbb{R}
            \end{array} $ \\

actions which, when extended to the solutions $U \in {\cal C}^\infty ( \Omega, \mathbb{R} )$
of the PDE in (A.1), will transform them into solutions of the same equation. \\

The well known difficulty here related to {\it global} actions of (A.2) on functions in ${\cal
C}^\infty ( \Omega, \mathbb{R} )$ is the following. \\
In general, the Lie group actions (A.2) defined on the Euclidean domains $M$ cannot so easily
be extended to act on the functions $U : \Omega \longrightarrow \mathbb{R}$ as well. And the
only problem here is that such extended actions {\it cannot} be defined so easily {\it
globally}, that is, for the functions $U : \Omega \longrightarrow \mathbb{R}$ considered on
the {\it whole} of their domain of definition $\Omega$, Rosinger [1, chapters 1,2]. The
reason for that is rather simple, namely, the lack of invertibility of certain functions
involved, Rosinger [1, pp. 14,15]. \\

For further clarity about the mentioned difficulty facing global Lie group actions on
functions, we recall here a simple example, namely, the {\it rotation in plane of a
parabola}. \\

Let us take the function $U : \Omega \longrightarrow \mathbb{R}$ given by

$$ U ( x ) ~=~ x^2,~~~ x \in \Omega = \mathbb{R} $$

and let us consider the Lie group $G$ on $M = \Omega \times \mathbb{R} = \mathbb{R}^2$ given
by the rotations of the plane around the origin $( 0, 0 ) \in \mathbb{R}^2$. Then it is
obvious that, unless it is an integer multiple of $\pi$, every such rotation, when applied to
all of the parabola, will turn it into a curve in plane which is {\it no longer} the graph of
any function in $V : \Omega \longrightarrow \mathbb{R}$. \\
Of course, bounded parts of the parabola can be rotated with sufficiently small angles, and
one again obtains the graph of a function. \\

In other words, arbitrary Lie group actions (A.2) {\it cannot} be extended to actions

$$ G \times {\cal C}^\infty ( \Omega, \mathbb{R} )
                      \longrightarrow {\cal C}^\infty ( \Omega, \mathbb{R} ) $$

This difficulty can, however, be easily overcome by the use of {\it parametric} representation
of the respective functions $U : \Omega \longrightarrow \mathbb{R}$, Rosinger [1, chapters
3-5]. In this way Lie group actions (A.2) can act {\it globally} on the functions $U : \Omega
\longrightarrow \mathbb{R}$ which are solutions of the rather general type of PDEs in (A.1). \\

For that purpose, we proceed as follows. Given any smooth function \\

(A.4) $~~~ U : \Omega \longrightarrow \mathbb{R} $ \\

we associate with it its {\it graph} \\

(A.5) $~~~ \gamma_U ~=~ \{~ ( x, U ( x ) ) ~~|~~ x \in \Omega ~\}
                                        \subseteq M = \Omega \times \mathbb{R} $ \\

Then by definition, a {\it parametric} representation of the function $U$ is given by any
smooth function \\

(A.6) $~~~ V : \Lambda \longrightarrow M $ \\

where $\Lambda \subseteq \mathbb{R}^n$ is nonvoid and open, such that \\

(A.7) $~~~ V ( \Lambda ) = \gamma_U $ \\

We note the following well known {\it advantage} of such parametric representations. Namely,
the set of functions in (A.6) is {\it larger} than that in (A.4). In other words, not every
function $V$ in (A.6) is the parametric representation of a function $U$ in (A.4). For
instance, a nontrivially rotated parabola in the plane can easily be written as a function in
(A.6), but not as a function in (A.4). \\

Therefore, we denote by \\

(A.8) $~~~ {\cal C}^\infty_n ( M ) $ \\

the set of all smooth functions in (A.6), and call them {\it n-dimensional parametric
representations in} $M$. \\

Clearly, we have the following embedding which associates with each function $U$ in (A.4) its
{\it canonical} parametric representation $V_U$ in (A.6), (A.7), namely \\

(A.9) $~~~ {\cal C}^\infty ( M, \mathbb{R} ) \ni U ~\longmapsto~
                                        V_U \in {\cal C}^\infty_n ( M ) $ \\

where for $U : \Omega \longrightarrow \mathbb{R} $, we define $V_U : \Omega \longrightarrow M$
by $V_U ( x ) = ( x, U ( x ) )$, with $x \in \Omega$. It follows that with the notation in
(A.6), we have in this particular case $\Lambda = \Omega$, therefore (A.7) holds, which means
that indeed $V_U \in {\cal C}^\infty_n ( M )$. \\

The important property of parametric representations is that for every Lie group $G$ acting on
$M$, see (A.2), one can naturally define the {\it global} Lie group actions on each of the
functions in ${\cal C}^\infty_n ( M )$, namely \\

(A.10) $~~~ G \times {\cal C}^\infty_n ( M ) \longrightarrow {\cal C}^\infty_n ( M ) $ \\

as follows. Given $g \in G$ and $V : \Lambda \longrightarrow M$ in ${\cal C}^\infty_n ( M )$,
we define \\

(A.11) $~~~ g V = g \circ V $ \\

where in the right hand term, $g$ denotes the mapping, see (A.2) \\

(A.12) $~~~ g : M \ni ( x, u ) \longmapsto g ( x, u ) \in M $ \\

while $\circ$ in (A.11) is the usual composition of mappings. \\

In other words, we define the action $g V$ in (A.11) by the commutative diagram \\

\begin{math}
\setlength{\unitlength}{0.2cm}
\thicklines
\begin{picture}(60,18)

\put(10,15){$\Lambda$}
\put(20,17){$V$}
\put(13,15.5){\vector(1,0){15}}
\put(30,15){$M$}
\put(40,17){$g$}
\put(34,15.5){\vector(1,0){15}}
\put(51,15){$M$}
\put(0,8){$(A.13)$}
\put(10.5,13){\line(0,-1){10.5}}
\put(10.5,2.6){\line(1,0){41.6}}
\put(52,2.6){\vector(0,1){10.5}}
\put(30,0){$g V$}

\end{picture}
\end{math} \\

And then, in view of (A.9), the simple construction in (A.13) allows the action (A.2) of every
Lie group on $M = \Omega \times \mathbb{R}$ to be extended {\it globally} to every smooth
function $U : \Omega \longrightarrow \mathbb{R}$. \\

{\bf Remark A} \\

From the point of view of {\it genuine Lie semigroups}, the essential feature of the
definition of {\it action on functions} in (A.11), (A.13) is that it is valid not only for Lie
group elements $g$, which therefore generate bijections (A.12), thus elements of  ${\cal
D}iff^\infty ( M, M )$. \\
Indeed, (A.11), (A.13) make also sense for {\it all smooth} mappings in the far larger ${\cal
C}^\infty ( M, M )$, thus for mappings generated by $g$ which need no longer be elements of
Lie groups, and instead can belong to genuine Lie semigroups as well.

\hfill $\Box$ \\

The essence of the above definition (A.10) - (A.13) of {\it global action on functions} is
very simple when seen in {\it categorial} terms, that is, in terms of most general properties
of the usual {\it composition} of functions. Indeed, initially, the functions of interest on
which the actions are supposed to be defined are, see (A.4) \\

(A.14) $~~~ U : \Omega ~\longrightarrow~ \mathbb{R} $ \\

while the actions operate according to, see (A.2), (A.3) \\

(A.15) $~~~ M ~\stackrel{g}\longrightarrow~ M,~~~ \mbox{with}~~~ g \in G $ \\

where $M = \Omega \times \mathbb{R}$. In this way, the extension of the actions (A.15) to
functions (A.14) leads to having to deal with the inversion of certain functions which may
fail to exist, Rosinger [1, pp. 14,15]. \\

However, if the functions (A.14) are embedded into the {\it larger} set of functions, see
(A.6) \\

(A.16) $~~~ V : \Lambda \longrightarrow M $ \\

where $\Lambda \subseteq \mathbb{R}^n$ is nonvoid and open, and this embedding is done
according to, see (A.9) \\

(A.17) $~~~ {\cal C}^\infty ( M, \mathbb{R} ) \ni U ~\longmapsto~
                                        V_U \in {\cal C}^\infty_n ( M ) $ \\

then the mappings (A.15) and (A.16) can trivially be composed with one another, thus yielding
(A.11), (A.13). And obviously, such a composition of mappings does {\it not} require the
mappings $g$ in (A.15) to be {\it bijections}, that is, to belong to ${\cal D}iff^\infty ( M,
M )$. Instead, they can belong to the far larger ${\cal C}^\infty ( M, M )$. \\

It appears that the above definition (A.10) - (A.13) of a {\it global} action on all smooth
functions by arbitrary Lie groups was presented for the first time in Rosinger [1, chapters
1-5], based on the above simple device of {\it parametric} representation of functions. \\ \\

\end{document}